\documentclass[11pt]{article}
\usepackage{amsmath,amssymb,amsthm}
\usepackage[greek,english]{babel}
\usepackage{graphicx}
\usepackage[round]{natbib}
\usepackage[continuous]{pagenote}
\makepagenote

\renewcommand*{\pagenotesubhead}[1]{}

\theoremstyle{definition}
\newtheorem{definition}{Definition}

\def\beq{\begin{equation}}
\def\eeq{\end{equation}}

\def\ter{{\rm\bf ter\,}}

\defcitealias{CP}{CP}
\defcitealias{NEM}{NEM}

\makeatletter
\date{}
\makeatother

\begin{document}

\title{Is Peirce's reduction thesis gerrymandered?}

\author{Sergiy Koshkin\\
\\
Department of Mathematics and Statistics\\
University of Houston-Downtown\\
1 Main Street, Houston, TX 77002\\
e-mail: koshkins@uhd.edu\\}
\maketitle
\begin{abstract} We argue that traditional formulations of the reduction thesis that tie it to privileged relational operations do not suffice for Peirce's justification of the categories, and invite the charge of gerrymandering to make it come out as true. We then develop a more robust invariant formulation of the thesis by explicating the use of triads in any relational operations, which is immune to that charge. The explication also allows us to track how Thirdness enters the structure of higher order relations, and even propose a numerical measure of it. Our analysis reveals new conceptual phenomena when negation or disjunction are used to compound relations.

\bigskip

\end{abstract}

\section*{Introduction}

Peirce's reduction thesis, henceforth PRT, is deeply rooted in the system of categories he developed after getting acquainted with de Morgan's logic of relations c.\,1867, and defended for the rest of his life.  The more contentious core of PRT can be stated as a conjunction of two clauses:
\begin{itemize}
\item {\bf Triadic irreducibility:} some triadic relations cannot be reduced to combinations of monadic and dyadic ones.

\item {\bf Polyadic reducibility:} all tetradic and higher adicity relations can be reduced to combinations of triadic ones.
\end{itemize}
\noindent After Peirce's death, aside from Quine's brief passage in the 1935 review of Peirce's Collected Papers, there were no direct references to PRT until 1970s. However, there did appear several counterclaims, notably in \citep{Low15} and \citep{Quine54}\pagenote{Quine only discussed the reducibility to triads clause in 1935, and made no mention of Peirce, PRT or what his 1954 reduction to dyads meant for it then, or later. Not in 1981, when he simplified the reduction procedure in the 4th edition of Mathematical Logic, and not in the 1995 essay Peirce's Logic. His silence is usually interpreted as a courteous submission that PRT is false, but he never did say so himself. The pairing construction was pressed as an objection to PRT by others, e.g. \citep{ChrJoh81}.}, to the effect that all polyads are reducible to dyads, as Kempe already suggested in Peirce's lifetime. The proposed reduction was different from Kempe's hypostatic abstraction, and was based on the now familiar set-theoretic representations of a triple as a pair of a unit and a pair, of a quadruple as a pair of pairs, and so on. Whether these reductions do what they appear to do turns on precisely what ``combinations" one allows in a  ``reduction", at least on the now classical formalizations of PRT by \cite{Herz} and \cite{Bur91}.

This turn of events was anticipated in what seems to be the first scholarly work on PRT specifically, \citep{Skid71}, where the term ``Peirce's thesis" is coined ("reduction" was added by Herzberger). After reviewing Kempe style reduction to dyads, Skidmore suggested that one could save PRT by restricting operations on relations, but that such a move would make it into an artifice of one's favored operations. 

Skidmore's concern, under the catchy name of ``gerrymander", recently resurfaced in \citep{Conr20}, no longer as a potential counter to hypothetical defenses, but as a pointed objection to Herzberger's and Burch's versions of PRT. Conarroe, who, unlike Skidmore, is sympathetic to the thesis, argues that a response equal to Peirce's ambition for the categories requires a meta-argument demonstrating that Thirdness is unavoidable no matter what base operations one adopts, and suggests Hintikka's formalization of Peirce's corollarial/theorematic distinction as a template. 

We aim to offer such a meta-argument based on a strengthened formulation of PRT. After arguing that the gerrymandering charge does have a bite against the humbler versions of the thesis, we allow {\it any} (first order) relational operations in a reduction, but {\it explicate} their use in a precise sense (Section \ref{PRT+}). That means that devices that merely abbreviate the use of triads are explicitly presented as doing so. In the standard first order logic (FOL) they are identifications of free and bound variables in different places\pagenote{That identifications conceal relations was suggested already by \cite{DeM60}:``In the consideration of the proposition, {\it identification} of {\it objects} is in truth a {\it relation} of {\it concepts}" (emphasis his).}, and in EG they are even more visible as branch points. The resulting invariant version of PRT is implicit in the state of the art mathematical proofs \citep{CorDau06,CorPos04,CorPos06,CorPos11}, and is not only stronger, but also more defensible against the gerrymandering and other objections. It gives truth to Peirce's unqualified conviction that ``triadic relation cannot be reduced without  a use of triadic relation" (\citetalias{CP}\,3.423-4).

It is ironic that, according to \cite{Mer79}, Peirce's argument for PRT ``turns on loose verbal reasoning that upon examination will not hold up", while his and many other arguments against it fail exactly because loose verbiage tends to conceal habitual uses of triads. Moreover, the explication shows that even formal languages, which we task with unearthing hidden moves and assumptions, are not entirely up to the task when it comes to PRT. For example, the pairing construction, as formalized in set theory, relies on the very FOL devices that are explicated as using triads. The same syntactic absorption is traceable in other formalisms, e.g. in the term-functor logic \citep{Som00}. As a result, proposed reductions to dyads, when taken in anything other than a narrow technical sense, emerge as instances of the use/mention error pointed out by Peirce in connection with dyadic algebra, ``the very triadic relations which it does not recognize it itself employs" (\citetalias{CP}\,8.331).

The explication also allows us to sketch how relational Thirdness manifests semiotically as informational dependence between the relata, and detect the differences in its operation in positive (with $\exists,\land$ only) and non-positive reductions. In a parallel to Hintikka's treatment of theorematicity, we introduce a numerical measure of genuine triadicity in a relation that we call {\it ternarity}, the minimal number of triads needed to represent it when their use is not concealed by syntactic absorption. We then describe {\it genuine} triads and polyads as those having maximal possible ternarity, and briefly touch on new directions of their study and classification. 

On the side of mathematics the case of PRT has been settled, and it has been settled in Peirce's favor. The work on PRT in the last two decades 
poses novel and challenging questions, 
and we hope to reframe PRT from a controversial perennial problem to an established foundation for philosophically (and mathematically) fertile research program in Peircean analysis of relations and categories.

\section{Gerrymandering objection}

Peirce first stated PRT\pagenote{Part of the thesis was vaguely anticipated already by \cite{DeM60}:``Two thoughts cannot be brought together in thought except by a thought: which last thought contains their relation". This mini-argument, that bringing two together presupposes a third, is often used by Peirce and later authors in defense of PRT.}, without much argument, in Description of a Notation for the Logic of Relatives (1870, \citetalias{CP}\pagenote{Standard abbreviations are used for Peirce's works: CP\,v.p -- The Collected Papers of Charles Sanders Peirce; NEM v:p -- The New Elements of Mathematics by Charles S. Peirce; v volume, p page.}\,3.144), and only gave it significant elaboration and defense after reflecting on Kempe's Memoir on the Theory of Mathematical Form (1886) \citep{Anel97}. A proof is sketched informally in A Guess at the Riddle (1890, \citetalias{CP}\,1.363), where we find the valency argument for the irreducibility and the use of hypostatic abstraction for the reducibility. It is refined in The Critic of Argument (1892, \citetalias{CP}\,3.423-4), where Peirce responds to Kempe's construction that seems to show reducibility to dyads. He reconsidered this response after developing the system of Existential Graphs (EG) since 1896, and we find Peirce's mature position in the Lowell lectures of 1903 (\citetalias{CP}\,1.345-7), and in a 1905 letter to James (\citetalias{NEM}\,3:832-3) prompted by Royce's objections. 

In A Guess at the Riddle the tetrad ``$A$ sells $C$ to $B$ for the price $D$" is analyzed into ``$A$ makes with $C$ a certain transaction $E$" and ``$E$ is a sale of $B$ for the price $D$". This $E$ is a newly minted {\it ens rationis} created by hypostatic abstraction. Formally, we wish to reduce $S(A,B,C,D)$\pagenote{Throughout the paper we abuse terminology and notation by ignoring the difference between relations as sets of tuples, and predicates that are interpreted as their truth conditions on a domain. As \cite{Kerr92} noted, this difference does not come into play in the context of reduction. In particular, we will use ``relation" and ``predicate" interchangeably.} to $S'(A,C,E)\land S''(E,B,D)$, but the latter has $E$ in it, which is not in the original $S$. What we really want to say is that there {\it is} such a transaction, so the reduction is
\beq\label{sell}
S(A,B,C,D)\equiv\exists x\left[S'(A,C,x)\land S''(x,B,D)\right].
\eeq
Peirce then explains why a similar trick would not work for reducing triads to dyads:``But when we attempt to imitate this proceeding with dual relatives, and combine two of them by means of an $X$, we find we only have two blank places in the combination" (\citetalias{CP}\,1.363). This is a rudimentary form of the valency argument fully developed in \citep{Herz}. As applied to the example from The Critic of Argument, ``$A$ gives $B$ to $C$", it means that we cannot similarly reduce  $G(A,B,C)$ to dyads because in $G'(\cdot\,,x)\land G''(x,\cdot)$, with any $G',G''$, there are only two dots, but there are three of $A,B,C$.

So far, so good. Except, as Peirce deduced from Kempe, we can do the following:
\beq\label{give}
G(A,B,C)\equiv\exists x\left[G'(x,A)\land G''(x,B)\land G'''(x,C)\right],
\eeq
where $x$ is a ``certain act", and $G':=$ ``in $x$ something is given by $A$", $G'':=$ ``in $x$ someone is given $B$", and $G''':=$ ``in $x$ something is given to $C$". To this, Peirce counters that this reduction ``fails to afford any formal representation of the manner in which this abstract idea is derived from the concrete ideas" (\citetalias{CP}\,3.424). The implication is, presumably, that if such a representation were given it would reveal hidden triadicity. 

This objection is often pressed against pairing reductions {\it a la} L\"owenheim and Quine as well, but it is a double edged sword, as Skidmore pointed out. Trouble is, what is good for the goose is good for the gander. If we are to fault \eqref{give} for the sin of omission we must fault \eqref{sell} for that same sin. And if \eqref{sell} is judged to conceal tetradicity then no reduction to triads is accomplished there either. We now see how improvident a position it was for Peirce to take -- accepting either horn of the dilemma, fault or no fault, unravels one of the clauses of PRT.

As Peirce realized by 1903, neither \eqref{sell} nor \eqref{give} fail to represent the relevant manners, and we rather need to drive a wedge between those manners as represented. However, the wedge Peirce drove in 1903 appeared to depend on {\it re}-representing the reductions in the calculus of Existential Graphs (EG) that he developed in the meantime. And this is the feature Skidmore seized upon:
{\small \begin{quote}Peirce's defense of his Thesis here amounts to showing that the system of existential graphs is fundamentally superior to the system of quantification theory; otherwise its only advantage might be construed as appearing to preserve Peirce's Thesis.
\end{quote} }
\noindent This is the gerrymandering objection. It is reprised and so nicknamed in \citep{Conr20} after reflecting on Herzberger's and Burch's proofs of PRT that adopted the response of Lowell lectures:``One can simply argue that the logical system itself -- i.e. EG -- is gerrymandered in a way to ensure something like the Reduction Thesis holds". Indeed, as Herzberger notes, the EG privilege certain operations on relations over others, namely the ``bonding" of predicates two at a time as in \eqref{sell}, and disallow the ``triple junction" in \eqref{give}. \cite{Bur21} similarly states that ``both Peirce and Quine were correct: the issue entirely depends on exactly what constructive resources are to be allowed to be used in building relations out of other relations".

However, such a humbled version of PRT would not suffice for Peirce's purposes of carving nature at the joints by his categories, with far reaching epistemic and metaphysical consequences. What would become of it if Quine can choose operations that reduce relations to dyads, someone else have them reduce to $n$-ads for some $n>3$, or even have irreducible relations of every finite adicity, as \cite{Anel97} reads Tarski to suggest?
It cannot simply be that Peirce and Quine are both right. For Peirce to be right Quine has to be {\it wrong} in a deeper conceptual sense, despite being correct in a technical algebraic sense.

To this end, we need an argument that triadicity, or rather Thirdness, will manifest in some form one way or another no matter what operations one chooses. Such an argument can be developed from an {\it invariant} formulation of PRT that does not depend, in particular, on a choice of base operations. It can be done, for example, by identifying Thirdness not only in relations, but in operations as well, or by identifying it in relations without a recourse to operations. It needs to be done if only to make sure that one is not allowed to cheat by simply folding the use of triads into operations, and it is attempted already by Herzberger and Burch. However, as we will argue, their ``Thirdness" in operations has a tenuous connection to what one might independently describe as such, and a very clear connection to what they need to prove PRT. This only reinforces the gerrymandering concern instead of dispelling it.

\section{Valency rule}

After decades of obscurity, counterclaims and objections \cite{Herz} presented a formalization, the bonding algebra, where both clauses of PRT could be proved. It was, in many ways, a toy model, but it was a proof of concept. Before him, not only was PRT held quaint or artificial, it was even doubted that its two clauses can both hold in a consistent framework \citep{Mer79}.

\begin{figure}[!ht]
\begin{centering}
a)\ \ \ \includegraphics[width=0.20\textwidth]{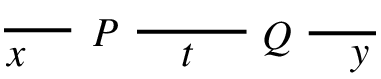} \hspace{0.1in} b)\ \ \ \includegraphics[width=0.10\textwidth]{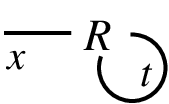} \hspace{0.1in} c)\ \ \ \includegraphics[width=0.21\textwidth]{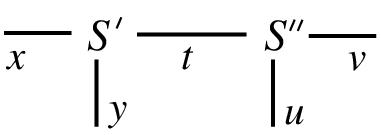} \hspace{0.1in} d)\ \ \ \includegraphics[width=0.13\textwidth]{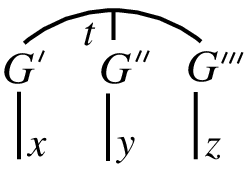}\par
\end{centering}
\caption{\label{EGex}Existential graphs for: a)\, bonding dyads $\exists\,t\left[P(x,t)\land Q(t,y)\right]$; b)\,self-bonding a triad $\exists\,t\,R(x,t,t)$; c)\,bonding triads $\exists\,t\left[S'(x,y,t)\land S''(t,u,v)\right]$;\ \ \ \ \ d)\, triple junction $\exists\,t\left[G'(x,t)\land G''(y,t)\land G'''(z,t)\right]$.}
\end{figure}

The bonding algebra is a simplified algebraic variant of EG, and for our purposes a simple fragment of them will suffice, corresponding to a fragment of first order logic (FOL) with only $\exists$ and $\land$. The relations are represented by predicate letters scribed on a two-dimensional {\it sheet of assertion} with lines radiating from them to represent their variable places\pagenote{Strictly speaking, Peirce used "existential graphs" to express closed sentences only, those expressing  predicates he called "valental graphs" \citep{Bur11}, but we will keep the more familiar label. To represent non-symmetric relations the lines should be numbered. One can place explicit labels on them, or agree, say, that the first place corresponds to the horizontal line to the left, and the rest go counterclockwise in order. Bookkeeping operations of permuting predicate places are included into the bonding algebra to reflect the graphs properly. Since the numbering and permutations will be immaterial to our discussion we will gloss over these devices.}. We will agree that loose lines on a predicate represent free variables\pagenote{This is in contrast to Peirce's use of loose ends for existential quantification. Single place quantification is not allowed in the bonding algebra anyway.}, and often write variable names next to them for convenience of translation between EG and FOL. Placing several predicates on the sheet amounts to taking their conjunction, and connecting two loose ends amounts to existentially quantifying over a dummy variable placed into the corresponding places, Figure \ref{EGex}\,a),\,b),\,c). The latter is the bonding operation of Herzberger that Peirce calls relative product and Burch calls join. Note that bonding can occur only two at a time, so one can neither quantify over a single place, nor over more than two. Hence not even every $\exists,\land$ FOL formula is representable, for example, the triple junction on Figure \ref{EGex}\,d) is not. Moreover, the final graphs must be connected.

The bonding algebra allows reducing $n$-ads to triads by iterating the bonding on Figure \ref{EGex}\,c), while blocking reduction to dyads. The latter is demonstrated by extending Peirce's valency argument. If two predicates have adicities (valencies) $\mu$ and $\nu$ then their bond has adicity $\kappa=\mu+\nu-2\lambda$ with $\lambda\geq1$, where $\lambda$ is the number of bonds. This is the {\it valency rule}. When $\mu,\nu\leq2$ obviously $\kappa\leq2$. No matter how many monads and dyads we bring in, and how many times we bond them, we will never get anything other than monads and dyads in the end. Thus, both clauses of PRT are proved\pagenote{We glossed over permutations and negation. Both of them preserve adicity, and hence obey the valency rule.}. 

Herzberger himself admits that this is ``too neat to be true". And so it is. Why is the triple junction not allowed? Because it breaks the valency rule. What of quantifying over a single place? Same reason. Even mere juxtaposition without bonding is not allowed, which is why we cannot get disconnected graphs. The valency rule is elevated to Thirdness itself, mathematically embodied, in Herzberger's account. However, while one may reasonably suspect triple junction of foul play, single place quantification or juxtaposition do not exactly scream Thirdness or triadicity. The valency rule seems like too broad a brush that paints Thirdness into corners with no tangible connection to it. 

Also, completely missing from Herzberger's picture is  teridentity, which played a major role in Peirce's later work on reduction. Something like the bonding algebra might have been what Peirce had in mind in 1890, but he had misgivings about it already in 1892, and moved past it by 1903.  

\section{PAL and teridentity}\label{PAL}

While Herzberger's paper inspired the more positive treatment of PRT in the 1980-s, the waterhsed came with the publication of \citep{Bur91}. Burch approached Peirce's work on relations and their reduction much more comprehensively, both in chronological and philosophical scope, and interpreted it in terms of modern algebraic logic. In particular, he defined syntax and semantics of a formal language PAL (for Peirce's Algebraic Logic) that translates EG into a more conventional form. \cite{Kerr92} and \cite{Bur97} give helpful informal expositions.

\begin{figure}[!ht]
\begin{centering}
a)\ \ \ \includegraphics[width=0.13\textwidth]{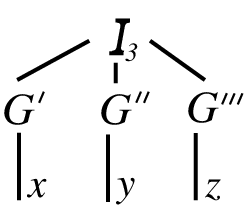} \hspace{0.1in} b)\ \ \ \includegraphics[width=0.10\textwidth]{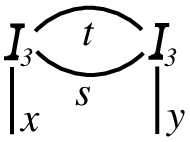} \hspace{0.1in} c)\ \ \ \includegraphics[width=0.31\textwidth]{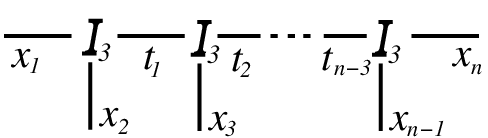} \par
\end{centering}
\caption{\label{I3show}PAL simulations via teridentity for: a) triple junction; b) binidentity; c) $n$-identity.}
\end{figure}

\begin{figure}[!ht]
\begin{centering}
a)\ \ \ \includegraphics[width=0.19\textwidth]{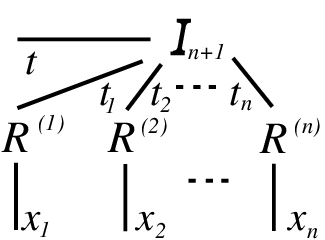} \hspace{0.05\textwidth} b)\ \ \ \includegraphics[width=0.50\textwidth]{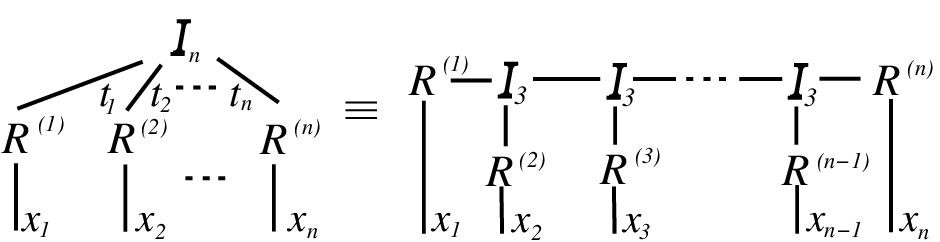} \par
\end{centering}
\caption{\label{Inshow} PAL simulations for a) identification of $n$ free variables; b) $n$-junction in hypostatic abstraction split into teridentities.}
\end{figure}

In addition to operations of the bonding algebra, PAL allows juxtaposition of predicates (Cartesian product of relations), albeit with a fateful restriction to be discussed later, and so it is no longer the case that monads and dyads cannot generate any triads at all. And, in agreement with late Peirce, it gives the pride of place to the teridentity predicate $I_3$, interpreted as a triad that contains all triples of identical elements, and them only. While the bonding algebra simply excluded triple junctions, PAL simulates them by using $I_3$, see Figure \ref{I3show}\,a). However, this no longer threatens PRT because the simulation explicitly involves a triad, namely $I_3$. The ``natural" definition of $I_3$ in terms of binidentity  
$I_3(x,y,z):=I_2(x,y)\land I_2(y,z)$, or in a more familiar form $(x=y)\land (y=z)$, is blocked in PAL, because it requires identification of variables in different places without quantification. In contrast, not only $I_2$ but pluridentities $I_n$ of any finite adicity $n$ are definable from $I_3$ by mere bonding, see Figure \ref{I3show} b), c).
With them at hand, one can simulate identification of variables in any number of places of any predicates, and junctions of any multiplicity, Figure \ref{Inshow}. 
\begin{multline}\label{multiidn}
\exists t_1\dots\exists t_n\left[I_{n+1}(t,t_1,\dots,t_n)\land R^{(1)}(t_1,x_1)\land\dots\land R^{(n)}(t_n,x_n)\right]\\
\equiv R^{(1)}(t,x_1)\land\dots\land R^{(n)}(t,x_n).
\end{multline}

Teridentity is clearly the power engine of PAL. One reason for its fertility is a feature it shares with all triads, noted already by Peirce. When a triad is bonded with other predicates the total number of loose ends (free variables) increases by $1$, and one can iterate to increase it indefinitely. In contrast, bonding with $I_2$ leaves adicity unchanged, and so binidentity is adicity sterile. However, the common terminology contributes to the habit of overlooking triadicity. Graphs, like the one on Figure \ref{I3show}\,c), are traditionally called {\it binary} trees despite the fact that their internal vertices are {\it trivalent}. This is closely related to how the pairing construction conceals triadicity as well.

The hypostatic abstraction can be performed in PAL without introducing any relation specific triads as in \eqref{sell}, teridentity is the only triad needed. To do so, given a relation $R(x_1,\dots,x_n)$ to be reduced, we add all $n$-tuples of elements in the domain to it as new elements, and define dyads $R^{(i)}(t,x)$ that contain tuple-elements from $R$ in the $t$ place and their $i$-th entries in the $x$ place. Then $R$ is reduced, ostensibly to dyads, as follows, see Figure \ref{Inshow} b):
\begin{multline}\label{hypred}
R(x_1,\dots,x_n)\equiv\exists t\left[R^{(1)}(t,x_1)\land\dots\land R^{(n)}(t,x_n)\right]\\
\equiv\exists t_1\dots\exists t_n\left[I_{n}(t_1,\dots,t_n)\land R^{(1)}(t_1,x_1)\land\dots\land R^{(n)}(t_n,x_n)\right],
\end{multline}
where we used that $\exists t\,I_{n+1}(t,t_1,\dots,t_n)=I_{n}(t_1,\dots,t_n)$. The second line shows, from the perspective of PAL, how $n$-identity, and hence teridentity, is surreptitiously employed in the traditional ``reduction to dyads". 

The set-theoretic pairing construction employs the same trick. To treat a triple $(x_1,x_2,x_3)$ as a pair $((x_1,x_2),x_3)$ we need to turn the pair into an element $t=(x_1,x_2)$, and relate its entries to that element, say, by dyads $R^{(i)}(t,x)$. We then need to relate the pair-element to the third entry, say, by another dyad $\widehat{R}(t,x)$. The reduction then amounts to the triple junction
\beq\label{pairred}
\exists t\left[R^{(1)}(t,x_1)\land R^{(2)}(t,x_2)\land\widehat{R}(t,x_3)\right]
\eeq
that asserts existence of such a pair-element, and conceals triadicity in the same way as \eqref{hypred}.

Burch's is an impressive account of polyadic reducibility and the singular role of teridentity in it, but what about irreducibility? Here comes the restriction on juxtaposition that Burch enforces by his apparatus of arrays and assemblies -- one can only juxtapose subgraphs at the very end. Since this is the only operation of PAL that breaks the valency rule reducible relations necessarily factor into Cartesian products of lower adicity relations. Following Peirce's terminology, Burch calls such relations {\it degenerate}, as opposed to genuine ones, like pluridentities. For triads, reducibility to monads and dyads is even equivalent to degeneracy. 

However, one cannot shake off the suspicion that operations of PAL are lined up just so to make it happen. Peirce, it is true, second graded juxtaposition/Cartesian product, but that was because he did not consider it a ``genuine" operation that combines concepts non-trivially, but rather more of a listing device \citep{Brun97}. If anything, this is a reason for including it at full strength when probing irreducibility, not excluding or restricting it, for a mere listing device is surely Thirdness-free. On the other hand, Burch's proof of irreducibility would not go through without such a restriction. It is certainly an improvement on the full-blown valency rule, but the shadow of gerrymander still looms over it.

\section{Discerning Thirdness}

To get a firmer grasp on what should and should not be allowed in a reduction, we need a clearer idea of how Thirdness manifests in the logic of relations. There is a wrinkle already with relations themselves. Not just any triad has Thirdness, says Peirce, only {\it genuine} triads do. Which triads are genuine? Peirce is somewhat ambiguous on the matter. Sometimes degeneracy means ``mere coexistence" (\citetalias{CP}\,1.372), i.e. juxtaposition, as in Burch. Sometimes, especially for triads, it means reducibility, which, of course, depends on the means of reduction.

The situation is even murkier with those means, relational operations. We know that Peirce privileged relative product (Herzberger's bonding), and sought to eliminate non-relative Cartesian and Boolean product (conjunction with identification of variables) from his algebras \citep{Brun97}. Later, triple and higher junctions joined the blacklist, but, with Cartesian product, at least, the reason for the blacklisting was triviality, not Thirdness.

Herzberger's valency rule was a more principled Thirdness criterion, but it swept too broadly. Burch proposed a different, PAL based, criterion:``Let us say that Thirdness is involved in any operation or procedure... such that, if the operation or procedure were formalized in PAL, its definition would have to presuppose the availability of at least one non-degenerate triadic relation" \citep[p.117]{Bur91}. In other words, an operation ``involves Thirdness" when one needs triads to express it in PAL.

This would be a good criterion if we trusted that everything that does {\it not} involve Thirdness is covered by PAL {\it sans} teridentity. Here the two clauses of PRT pull in the opposite directions. The weaker the means we need to reduce polyads to triads the stronger the claim for locating the reducing power in triads themselves. However, this reverses for  irreducibility. That triads are irreducible by PAL's meager means is not so impressive, especially when the means left out are not particularly triadic. One can plausibly attribute irreducibility to the poverty of means rather than to special powers of triads, contrary to Peirce's intention in PRT. And if we suspect that PAL operations are handpicked to highlight teridentity then detecting Thirdness by Burch's criterion is circular -- operations that were left out will ``involve Thirdness" by design. 

For example, according to Burch's criterion, any operation that produces odd-adic relations from even-adic ones exclusively ``involves Thirdness". In particular, this implicates quantification over a single variable. However, we can express it by bonding with self-identity monad $I_1$ (true on every element in the domain): $\exists x\,P(x,y)\equiv\exists x\left[I_1(x)\land P(x,y)\right]$. True, one needs $I_3$ to express $I_1$ in PAL, by bonding two of its loose ends as on Figure \ref{EGex}\,b), but teridentity itself can be similarly expressed by bonding two ends of pentidentity. That does not mean that it involves Fifthness.

To summarize, both Herzberger's and Burch's criteria are too liberal in attributing Thirdness, but there is a natural PAL based criterion that is more conservative. If we suspect that PAL$-I_3$ is too weak itself then an operation that allows to reduce teridentity when added to it {\it must} carry Thirdness within. On this criterion, single place quantification {\it does not} involve Thirdness. Indeed, if adding it to PAL$-I_3$ could reduce teridentity to monads and dyads then so could PAL$-I_3$ by itself. After all, we can use any monad in a reduction, including $I_1$. This already gives us a strengthened version of PRT's irreducibility clause: triads remain irreducible even if we add single place quantification to PAL. On the other hand, variable identification in distinct places is a different matter. Since $I_3(x,y,z):=I_2(x,y)\land I_2(y,z)$ it is laced with Thirdness. The same goes for triple junction, as we already saw.

Thus, our conservative criterion seems to gauge Thirdness in operations better, but it is, at best, an indicator. What we still lack is a handle on what it indicates, exactly. Moreover, it raises another concern. Bonding, as in $\exists xR(x,y,x)$, looks like identification $R(x,y,x)$ followed by quantification $\exists x$. If identification is a source of Thirdness, and it is a step in bonding relations, then how come bonding itself remains Thirdness-free? If anything, it builds on it, it seems. Can it just be an artifact of FOL representation? 

Let us look at a different representation of relations that came up naturally in applications. \cite{Codd70} introduced the relational model of databases that deeply involves algebra of relations, albeit without being aware not only of Peirce-Schr\"oder's work, but even of Tarski's. It eventually entered Peircean literature after \cite{Sowa76} related database queries to conceptual graphs that are quite reminiscent of EG\pagenote{Sowa did not mention Peirce or EG in the 1976 paper, but described conceptual graphs as ``based on" EG in later works.}. In particular, more recent works on PRT, like \citep{CorDau06}, are informed by conceptions developed in the context of relational databases.

In the database model, relations are visualized as tables with columns representing predicate places, and rows containing the member tuples. Operations on relations are conceived as database queries that return relevant data organized into smaller table(s). They can be expressed either logically, in what Codd calls relational calculus, or procedurally, in relational algebra. One natural type of query deletes an entire column and prunes the duplicate tuples that are left. In relational calculus this is expressed by quantification, for example, $\exists xR(x,y,z)$ stands for deleting the first column of $R$'s table (a pair occurs in the output table whenever there was a triple in the original table with that pair in its second and third columns). 

Another type of query, called selection along two columns, returns rows with identical entries in them, keeping just one of the now duplicate columns \citep{Buss01}. The selection along the first and third columns of $R$'s table is expressed by $R(x,y,x)$. Thus, to bond along two columns we must, first, select rows where they have identical entries, and, second, delete the columns. We end up with the same identification followed by quantification procedure that we saw in FOL, but it is hard to charge the table model, unlike FOL, with being artificial. Nonetheless, Codd's relational algebra still splits PAL's bonding into more basic operations that involve Thirdness by any criterion. And $I_3$ is reducible in it. 

This is our first clue. The procedural view of FOL, made explicit in database queries, is clearly not the relevant one when it comes to Thirdness. In contrast, Peirce always emphasized {\it conceptual} combination and its dependence on what is combined. Identification of two variables establishes such dependence among {\it three} items, the two old variables and the new merged one, as \cite{DeM60} already noted. This triple dependence is then baked into the compound relation and has to be kept track of in any conceptual use of it. 

Not so with bonding. Even if we think of it as identification followed by quantification, the third item is quantified over and removed from consideration. There is no need to keep track of it anymore. Rather than building on the Thirdness of identification, quantification razes it to the ground. In the end, the bond retains dependence on only {\it two} items, the two old variables, and does not relate them to a third, even if, procedurally, we must introduce it momentarily to bond relations' extensions in practice. 

And thereby lies our second clue as well. Thirdness is not a conserved ``substance" that inheres in relational compounding. Being ``involved" is too vague a term, it can be increased (e.g. by identification), preserved (e.g. by bonding), or decreased (e.g. by quantification). Which case it is may also depend on context, although juxtaposition and bonding never increase it. It may even be possible to define ``amount of Thirdness" in a relational compound, and track how it is affected by operations. 

To summarize, Thirdness has much to do with conceptual dependence in relational compounds. It is affected by relations themselves, with their adicity as a (partial) indicator, but also by identification of predicate places and quantification in operations. This is the insight we will build on to devise a more comprehensive version of PRT to be fleshed out in the following sections.

\section{Invariant reduction thesis: positive irreducibility}\label{PRT+}

We will deal with positive  ($\exists,\land$) irreducibility first, because it is more transparent conceptually, and postpone including negation until Section \ref{neg}. Positive irreducibility is also interesting in its own right, after all, higher polyads are reducible by even weaker means. In contrast to Burch's PAL, quantification and juxtaposition will not be restricted in reductions, as in the upgraded version of PAL from  \citep{CorPos04}. The notion of explication introduced in this section will be central to all invariant formulations of PRT.

\begin{figure}[!ht]
\begin{centering}
a)\ \ \ \includegraphics[width=0.12\textwidth]{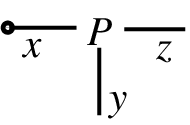} \hspace{0.1in} b)\ \ \ \includegraphics[width=0.24\textwidth]{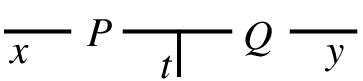} \hspace{0.1in} c)\ \ \ \includegraphics[width=0.40\textwidth]{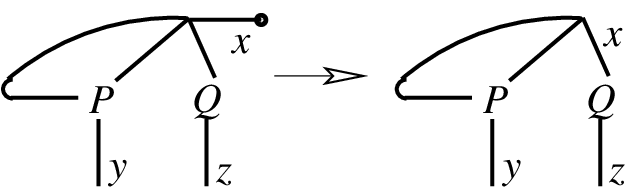}\par
\end{centering}
\caption{\label{EGmore} EG for a) single place quantification $\exists\,x P(x,y,z)$; b) free variable identification $P(x,t)\land Q(t,y)$; c) identification followed by quantification retracts to junction $\exists\,x\left[P(x,y,x)\land Q(z,x)\right]$.}
\end{figure}

Some EG devices that we left out so far will now be helpful for visualization. Single place quantification will be displayed by placing a heavy dot on the corresponding loose end\pagenote{As already mentioned, Peirce reserved loose ends for this purpose. \cite{Bur16} uses small circles, which reflects his idea that quantification ``involves" Thirdness in the form of teridentity with two bonded ends.}, Figure \ref{EGmore}\,a). When multiple lines connect to a single dot, all the corresponding places are quantified at once, as in triple junction. Identification of predicate places will be displayed by connecting their loose ends to a single branch point with a new loose end coming out of it\pagenote{As in \citep{Bur97}, where identification is a derivative operation of PAL defined via teridentity and called HOOKID.}, Figure \ref{EGmore}\,b). That a junction is identification followed by quantification means that a dotted end can be retracted to the branch point with the dot moved there, Figure \ref{EGmore}\,c). Since the branch points are already distinctive enough the dots are redundant, except at single ends, and we omit them in agreement with the conventional notation \citep{Rob73}.

\begin{figure}[!ht]
\begin{centering}
\ \ \ \includegraphics[width=0.60\textwidth]{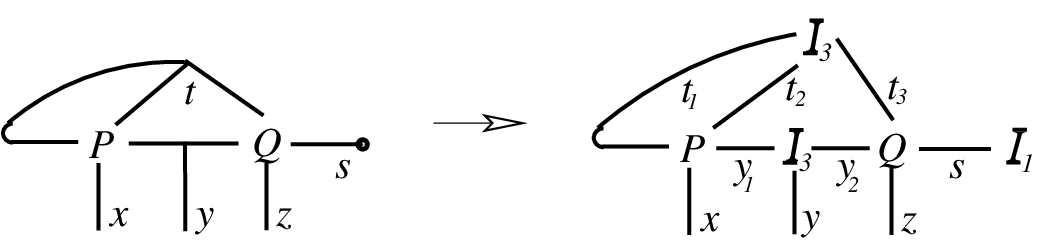}\par
\end{centering}
\caption{\label{EGexpli} Explication of dots and branch points in EG: $\exists\,t\left[P(t,x,y,t)\land\exists\,s\,Q(y,z,s,t)\right]$ to 
$\exists\,t_1\exists\,t_2\exists\,t_3\,\exists\,y_1\exists\,y_2\left[I_{3}(t_1,t_2,t_3)\land I_{3}(y_1,y,y_2)\land P(t_1,x,y_1,t_2)\land\exists\,s\left[I_1(s)\land Q(y_2,z,s,t_3)\right]\right]$}
\end{figure}

A key observation is that we get logically equivalent graphs by replacing dots and branch points by pluridentities of matching adicities, see Figure \ref{EGexpli}. This is, essentially, how the derivative operations that dots and branch points represent are defined in PAL\pagenote{Namely, QUANT and HOOKID in the notation of \citep{Bur91}.}. In other words, they are merely graphical aliases for bonding with pluridentity predicates. In particular, triple junctions and double identifications are stand-ins for teridentity, whether expressed in EG or in FOL. As Peirce put it, ``you may think that a node connecting three lines of identity is not a triadic idea. But analysis will show that it is so" (\citetalias{CP}\,1.346).   

And why stop at PAL operations? {\it Any} FOL formula or EG represents an operation on relations featured in it. These are called first order \citep{Buss01} or logical \citep{CorPos04} operations. 
\begin{definition}
We will call a logical operation {\it explicated} when the branch points and dots in its EG (or identifications and quantifications in the corresponding FOL formula) are replaced by pluridentities of matching adicity.
\end{definition}
\noindent If one prefers to explicate directly in FOL here is how. When a bound variable, say $t$, occurs in three or more different places, replace $t$ with a different $t_i$ in each place, parenthesize the entire scope of its quantifier, and replace $\exists\,t$ in front with $\exists t_1\dots\exists t_n\,I_{n}(t_1,\dots,t_n)\,\land$. When $t$ is a free variable occurring in {\it two} or more places, add $\exists t_1\dots\exists t_n\,I_{n+1}(t,t_1,\dots,t_n)\,\land$ in front of the entire formula after replacing $t$'s occurrences. The order in which this is done for different variables is immaterial.

One cannot say that explication privileges one notation over another because all it does is unfold devices that function exactly as certain predicates, namely pluridentities, in any notation. Surely, one is not allowed to discount the use of predicates in analysis on the `compelling' grounds that their use has been duly renamed into something else. That would be cheating. 

However, explication does clarify the significance of EG and PAL. They uniformize seemingly different means of relational compounding (in the positive fragment, for now) by reducing them to bonding predicates. With that, their respective contributions are made more visible, but they are operative within regardless of which operations or notation one uses to parse the compounding. We need not argue, even for broader philosophical purposes, that EG or base operations of PAL are superior in some “fundamental” sense, as Skidmore alleged, but rather only in {\it expository} sense. Similarly, eigenbases of linear operators are superior in that they express their operation in the most transparent form, but the same operation remains expressed in any other basis, it is an invariant. 
\begin{definition}
Let us call a relation {\it reducible} by a logical operation when it is represented by applying it to relations of strictly lower adicity. 
\end{definition}
\noindent Recall that a relation is degenerate when it is a Cartesian product of relations of strictly lower adicity, which is a special case of positive reduction. The irreducibility clause of PRT can now be formulated invariantly as follows. 
\begin{quote} {\bf Positive irreducibility (invariant form):} Non-degenerate triads are irreducible by any explicated positive logical operations.
\end{quote} 
Since degenerate triads are Cartesian products of monads and/or dyads this means that triads are positively reducible {\it if and only if} they are degenerate. The proof derives from a simple topological fact. Consider any representation of the triad by a (positive explicated) logical operation, and collapse all predicates in the representing EG into vertices. Place vertices at the loose ends as well. Those will be {\it pendants}, as vertices of valency $1$ are called  in graph theory. The end result is in an ordinary graph\pagenote{More precisely, it is called {\it pseudograph} in graph theory because it can have self-loops, when two loose ends of the same predicate are bonded, and multi-edges, when two predicates are linked by two or more bonding lines.} with at least three pendants, by assumption. Also by assumption, the graph is connected. If it was not then it would be a juxtaposition of two subgraphs that split the loose ends among them, and hence the relation would be a Cartesian product. 

The topological fact is that a connected graph with $3$ or more pendants must have a vertex of valency $3$ or more. It follows from Listing's census theorem oft-cited by Peirce\pagenote{To see this directly, consider paths going from one of the pendants to the other two. They share at least one edge, because there is only one coming out of the initial pendant, and they cannot share them all. Consider the last edge that they still share. At the end of it is a vertex at which the paths must diverge. Therefore, it has $1$ edge coming in and at least $2$ coming out, so its valency is $3$ or more.}. That trivalent vertex must come from a relation of adicity three or more. Thus, no representation of our triad is a reduction to monads and dyads only. 

Higher polyads can, of course, be reducible even when they are non-degenerate, but we can generalize PRT to them by considering only reductions to relations of adicity no more than three. Then a similar topological argument shows that for non-degenerate $n$-ads a positive reduction must involve at least $n-2$ triads. We will use this observation to propose a way of quantifying Thirdness in relations and in operations on them.

\section{A measure of Thirdness?}\label{meas}

In this section, based on the invariant form of PRT, we propose a numerical measure of relational Thirdness called positive ternarity, and use it to analyze what reductions to triads tell us about Peirce's  ``clustering together [of] ideas" in relations.

Before giving the definition, we need to address a subtlety with the reducibility clause of PRT. As presented in Section \ref{PAL}, it only guarantees that a higher polyad can be positively reduced to triads and dyads if the domain is extended. In other words, it is not the original relation that is reduced but some other (closely related to it, to be sure), on a larger domain. \cite{Herz} proposed a modification, where instead of adding tuples to the domain we put them in $1$-$1$ correspondence with some of the existing elements, and adjust the definition of $R^{(i)}(t,x)$ accordingly. Then the hypostatic abstraction variable runs over the original domain, and we have reduced the original relation. 

Alas, this only works for {\it small relations}, those with no more tuples than elements in the domain. For {\it large relations} we cannot be sure that a positive reduction to triads exists at all. Fortunately, any relation on an infinite domain is small\pagenote{Assuming the axiom of choice. In fact, that the cardinality of any infinite set is the same as the cardinality of its Cartesian powers is equivalent to the axiom of choice.}, and so are all pluridentities on any domain. Peirce did not seem to be concerned with extending domains whenever necessary, which effectively means that he thought of domains as indefinitely extended. We will restrict our discussion here to small relations and take up large ones in Section \ref{IPRTGen}.
\begin{quote} {\bf Positive reducibility:} Any small relation of adicity three or more is reducible to triads by explicated logical operations involving only $\exists,\land$.
\end{quote} 
Essentially, we want to count the number of triads in explicated EG that represent a relation. Of course, this number may vary from one EG representation to another, so  what we really want is the {\it minimal} number over all possible representations. 
\begin{definition}\label{ter+}
An EG is called {\it subtrivalent}\pagenote{We appropriate here a concept from graph theory, where a graph is called subtrivalent if all its vertices have valency at most $3$.} if all predicates in its explicated form have adicity at most three. {\it Ternarity} of a subtrivalent EG is the number of triads in its explicated form. {\it Positive ternarity}, $\ter^+$, of a relation is the minimal number of triads in its representation by a positive subtrivalent EG.
\end{definition}
\noindent The reducibility and irreducibility clauses of positive PRT state exactly that $\ter^+$ is well-defined and non-trivial, respectively. We emphasize that it can be strictly less than the number of triads in the representing EG. For example, deduction rules may transform it into an equivalent one with fewer triads, or new relations introduced through hypostatic abstraction may allow us to do the same. As a result, $\ter^+$ is not plainly readable off of EG, it is a {\it relation invariant}, like knot or graph invariants that are characteristics of knots and graphs independent of their representation.

This said, $\ter^+$ of monads and dyads is obviously $0$, and of teridentity or any other non-degenerate triad is $1$ (it cannot be more since any triad represents itself). The $n$-adic generalization we mentioned at the end of the last section means that positive ternarity of non-degenerate $n$-ads is at least $n-2$. In fact, it is exactly $n-2$, because hypostatic abstraction \eqref{hypred} unfolds any (small) $n$-ad into exactly $n-2$ teridentities that reduce $I_n$, see Figure \ref{Inshow}\,b). In particular, $\ter^+(I_n)=n-2$.

With a measure of Thirdness in relations at hand, we can better gauge it in operations as well, by looking at how they affect positive ternarity of relations operated on. There are many subtleties that we will only briefly mention here. Adding a triad or a triple point to EG does not necessarily increase its relation's ternarity. One may be able to eliminate it by consolidating a part of EG with it into a dyad or even a monad. Moreover, when there are many more triads than loose ends, hypostatic abstraction will reduce their number by completely replacing the original EG.

\begin{figure}[!ht]
\begin{centering}
a)\ \ \ \includegraphics[width=0.33\textwidth]{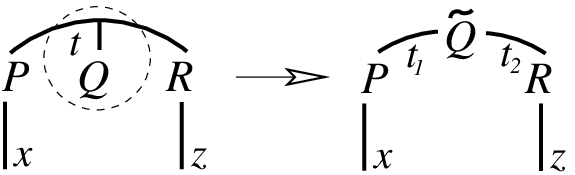} \hspace{0.5in} b)\ \ \ \includegraphics[width=0.29\textwidth]{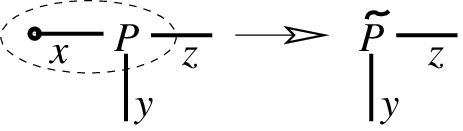} \par
\end{centering}
\caption{\label{Terop} a) Eliminating triple junction when one of the joined predicates is a monad by consolidating subgraph into a dyad, $\widetilde{Q}(t_1,t_2):=\exists\,t\left[I_{3}(t_1,t,t_2)\land Q(t)\right]$; b) ternarity reduces by $1$ when a free variable is quantified, $\widetilde{P}(y,z):=\exists\,xP(x,y,z)$. Dashed lines indicate consolidation and are not part of the EG.}
\end{figure}

As a result, the effect depends on the adicity of participating relations, and on whether or not they are degenerate. For example, triple junction will increase ternarity when applied to three non-degenerate dyads, but not if even one is degenerate or replaced by a monad, Figure \ref{Terop}\,a). It also depends on how the representing predicates are pre-bonded, or embedded into a larger EG. If we are looking for `pure' effect we should use 'generic' non-degenerate  predicates of sufficiently high adicity that are detached from a larger EG and each other. With that, triple junction and identification increase ternarity by $1$, bonding leaves it unchanged, and (single place) quantification decreases it by $1$, see Figure \ref{Terop}\,b).

Although $\ter^+$ may be a useful measure of Thirdness, it is a crude one even aside from accounting only for positive reductions. This takes us back to figuring out what relational Thirdness really is. On the verge of inventing EG, in On the Logic of Quantity (1896, MS\,13\pagenote{MS \# refers to the manuscript number \# in the Robin catalogue of Peirce's works.}), Peirce talked about it as  ``clustering together ideas into sets" that ``throws thoughts together" and manifests in reasoning through ``relational compounding". EG visualize how this clustering is mediated by non-degenerate triads that hold the compounds together, and suggest the idea of information flow.

We can think of loose ends as information channels that receive information when individuals (domain elements) are specified for them, fully or partially. If $R(x,y,z)$ is truly clustering and associating information then specifying $x=a$ should non-trivially influence other channels/places, by narrowing down which individuals can fill them under this relation\pagenote{In a slogan, ``information equals elimination of uncertainty" \citep[p.\,228]{Hin73}. And relation's Thirdness manifests in its ``welding", to use Peirce's word, of information flows. This extensional conception of information is too narrow on Peirce's view, but it gives us tentative means of quantifying it.}. In other words, the channels must be interdependent. This is reminiscent of Peirce's favorite semiotic triad of object-sign-interpretant, where the triad (sign) conveys the ``form" of one relatum (object) to another (interpretant). 

For pluridentities this interdependence is maximal, a perfect correspondence: specifying a single place determines all of them. For degenerate relations (Cartesian products) it is minimal, some channels are completely independent of each other. If $R(x,y,z)=P(x,y)\land Q(z)$ then specifying $x=a$ imposes no restrictions on $z$. This is the idea behind the proof that teridentity is non-degenerate on a domain with at least two elements $a\neq b$ \citep{Kerr92}. If the first two places can be filled independently of the third then $R(a,a,b)$ should be true if $R(a,a,a)$ and $R(b,b,b)$ are, and then $R$ is not teridentity for any $P,Q$. This reveals degeneracy as a semantic notion rather than algebraic one, like reducibility. Of course, there are intermediate grades between independence and perfect correspondence, but our $\ter^+$ is too crude to distinguish them. Any non-degenerate triad has as much ternarity as teridentity. 

In Shannon's information theory, dependence, or information sharing, between channels $X$ and $Y$ is measured by mutual information $I(X;Y)$, which is zero when they are independent and maximal when they determine each other. On finite domains, at least, we can assign equal probability to all tuples of a relation and then quantify mutual information between its different places and their sets. Taking the minimal value over all possible bipartitions of its places, the ``cruelest cut", will give us a finer measure of how close the relation comes to being degenerate. Such measures of ``information integration" are studied in biology \citep{Teg16}.

In the relational database theory functional and more general dependencies between columns (relation places) are studied since Codd to optimize data storage and facilitate access to it by losslessly decomposing relations into their projections. Such decompositions are closely related to Peircean reductions (the difference is that identified variables are not quantified over), and relies on exploiting dependencies between columns as well \citep[ch.\,7]{Maier83}. Clarifying how these dependencies are operative in EG representations corresponding to table decompositions should provide more clues on the nature of relational Thirdness.

\section{Negation and degeneracy}\label{neg}

As with the bonding algebra, the neat relationship between degeneracy and positive irreducibility is too neat to be true. Negation, or even already disjunction, complicate matters considerably. It is all the more glaring that no such complication is visible in Peirce's, or Herzberger's and Burch's treatments. Peirce's omission is probably explained by his focus on objections like Kempe's and Royce's, that questioned even {\it positive} irreducibility. And with Herzberger and Burch the explanation lies with the {\it other} operations they allowed in reductions. Negation does not cause problems in the bonding algebra because Herzberger's valency rule excludes other operations so sweepingly that negation has no partners to cause problems with. The valency rule rules out reduction of even degenerate triads. Burch's device is more subtle, but achieves the same end by allowing to negate connected EG only, not juxtapositions. As a result, degeneracy = reducibility holds for triads in Burch's system even without adding ``positive" to the latter. It is only in recent work, where such artificial restrictions are removed, that the proof of irreducibility became the ``difficult part" \citep{CorPos06,CorPos11}. 

The negation of even the most degenerate triad, like $\neg\left(P(x)\land Q(y)\land R(z)\right)$, is typically non-degenerate. By the de Morgan law, the same goes for free disjunctions (Cartesian sums) like $P(x)\lor Q(y)\lor R(z)$. A striking consequence of Burch's definitions is that they are thereby irreducible as well, despite what we see in the formulas with our lying eyes. The irreducibility that can be spoken thus is not the true irreducibility that Peirce intended. Burch softens the blow somewhat by proving that negated teridentity $\neg I_3$ is irreducible as well. So, at least, we can be sure that $I_3$ is not a Cartesian sum of monads and/or dyads. However, with a cleverer use of $\neg$ and $\lor$ all bets are off.

\begin{figure}[!ht]
\begin{centering}
a)\ \ \ \includegraphics[width=0.21\textwidth]{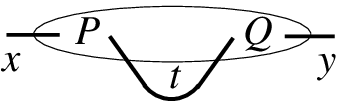} \hspace{0.0in} b)\ \ \ \includegraphics[width=0.13\textwidth]{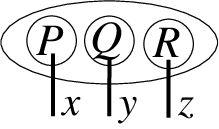} \hspace{0.0in} c)\ \ \ \includegraphics[width=0.42\textwidth]{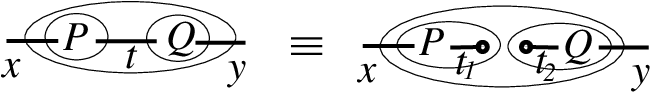} \par
\end{centering}
\caption{\label{EGcut} EG with cuts for a) bonding outside of negation $\exists\,t\,\neg\!\left[P(x,t)\land Q(t,y)\right]$; b)  the Cartesian sum of monads $P(x)\lor Q(y)\lor R(z)$; c) distributivity of $\exists$ over $\lor$, see \eqref{disE}.}
\end{figure}

Although $\exists,\land,\lor$ system is strictly weaker than logically complete $\exists,\land,\neg$, in EG $\lor$ and $\forall$ are expressed through negation anyway, by the device of cuts. The negated part of EG is enclosed by a simple closed thin line that bonding lines can cross but not end on. Cuts can nest, but not intersect each other. Sometimes it is necessary to pull branch points and/or bonding lines outside of a cut to reflect the correct order of quantification and negation, Figure \ref{EGcut}\,a). Loose ends are always extended outside the outermost cut. With the addition of cuts all restrictions on (first order) logical representation are removed. 

In the positive fragment we had a triple equivalence of degeneracy, disconnection and reducibility. This breaks down in the presence of negation, or even just disjunction. Cartesian sums have disconnected EG, but they are typically non-degenerate, Figure \ref{EGcut}\,b). At least, they are not {\it multiplicatively} (conjunctively) degenerate. If we wish to rule out {\it additive} (disjunctive) degeneracy for genuine triads as well we can make disconnection the basis of a new definition. 

Lest it be thought that using disconnection privileges EG, let us note that the concept is just as inherent in FOL. Let us call two variables linked when they share a predicate, and connected when they are the extremes of a chain of linked variables. Then the FOL formula obtained from the EG by Peirce's {\it endoporeutic} (outward in) reading method  \citep{Rob73} has disconnected variables whenever the EG was disconnected\pagenote{In the method, parts of EG outside of the outermost cut are conjoined with parenthesized expressions for enclosed parts that receive $\neg$ in front of them. The enclosed parts are read off the same way, iteratively. There are subtleties with reading the bonding lines when they cross cuts that we will not dwell on, but they link only at predicates or branch points. In the former case the corresponding variables share a predicate, and in the latter FOL uses the same variable in different places. There is one exception to this rule for so-called clipped lines that cross odd number of cuts entirely, but then the distinct variables will be linked by negated binidentities in FOL \citep{Shin00}. Thus, the resulting formula will be connected if and only if the EG was.}.

Of course, disconnected FOL formulas may be logically equivalent to connected ones, either by the artifice of introducing degenerate predicates, or, less trivially, by applying deduction rules like distributivity of $\exists$ over $\lor$, see Figure \ref{EGcut}\,c):
\beq\label{disE}
\exists t_1P(x,t_1)\lor\exists t_2\,Q(t_2,y)
\equiv\exists t\left[P(x,t)\lor Q(t,y)\right].
\eeq
However, disconnected FOL formulas can always be translated into disconnected EG, e.g. by Burch's procedure\pagenote{We first put the formula into the prenex form, taking care not to merge distinct variables, so that the prenex form of \eqref{disE}, for example, will be $\exists t_1\exists t_2\left[P(x,t_1)\lor Q(t_2,y)\right]$. Then we express $\forall$, $\lor$ and other connectives, if any, in terms of $\exists,\land,\neg$. The quantifier-free matrix is then scribed straightforwardly on the sheet of assertion. The lines coming out of the predicates are stretched out and linked to branch points placed within the corresponding cuts. If some variables were disconnected in the formula the corresponding lines will be disconnected in the EG, and explication will not alter that fact.} \citep{Bur97}. 
\begin{definition}
A relation is called {\it disconnected} if it is representable by a disconnected FOL formula or EG with free variables split among different connected clusters. It is {\it connected} otherwise.
\end{definition}
\noindent Disconnection generalizes both multiplicative and additive degeneracy. And it is a bridge to a version of general irreducibility clause of invariant PRT.

\section{Invariant reduction thesis: general irreducibility}\label{IPRTGen}

The same topological argument as before proves that reducible (by {\it any} explicated logical operation) triads are disconnected. The contrapositive gives us a version of invariant PRT. 
\begin{quote} {\bf General irreducibility (invariant form):} Connected triads are irreducible by any explicated logical operations.
\end{quote} 
This is not the difficult part. In the positive fragment, disconnection was degeneracy, and any positive reduction of triads produced a reduction of a very special form, the Cartesian product. Whether a relation is a Cartesian product was easy to check directly. What we still lack is a similar semantic description of disconnection in the presence of cuts, so it is hard to tell which relations are connected. Being representable by connected EG does not mean much, they need to be representable by connected EG {\it only}.

The $\lor\land$-decomposition introduced in \citep{CorDau06} and \citep{CorPos06,CorPos11}\pagenote{It is phrased there in set-theoretic terms and is, accordingly, named $\cup\cap$-decomposition.} provides just such a description. Disconnection partitions relation's places into non-overlapping clusters. While we may not be able to split its EG into a juxtaposition because of cuts, it turns out that it is a disjunction of finitely many relations whose EG do split into such a juxtaposition. 
\begin{definition}
A disconnected FOL formula is in {\it partitioned disjunctive form (PDF)} over a partition of its free variables into non-overlapping clusters if it is a disjunction of conjunctions with variables from a single cluster only in each conjunct. {\it Partitioned conjunctive form (PCF)} is defined similarly as a conjunction of disjunctions over a common partition.
\end{definition}
\noindent An example is 
\beq\label{oranddec}
R(x,y,z)\equiv\bigvee_iP_i(x,y)\land Q_i(z)
\eeq
with two clusters, $x,y$ and $z$. Using  distributivity laws, one can convert PDF into PCF and vice versa. Here is how it works for two disjuncts:
\begin{multline}\label{distrib}
\left(P_1(x,y)\land Q_1(z)\right)\lor\left(P_2(x,y)\land Q_2(z)\right)\\
\equiv\left(P_1(x,y)\lor P_2(x,y)\right)\land\left(Q_1(z)\lor Q_2(z)\right)\land\left(P_1(x,y)\lor Q_2(z)\right)
\land\left(P_2(x,y)\lor Q_1(z)\right).
\end{multline}
In the first two conjuncts variables from one of the clusters are missing, which we allow\pagenote{One can add fake predicates (that always return false) to them with the missing variables, if that is desired.}. As one can see from the formula, PDFs are always disconnected over the common partition. The result we need is that 
any disconnected FOL formula can be put into PDF (and PCF) over the same partition. The original proof explicitly appealed to iterative EG construction by PAL operations, which clashes with our invariant designs, but a direct argument in FOL can be given instead\pagenote{Put the formula into a prenex form while preserving the disconnection. If there are no quantifiers in the prefix, or the innermost quantifier is $\exists$, put the matrix into DNF, and in each disjunct group together predicates with variables from the same cluster. This will go through because  variables from different clusters never share a predicate. Then distribute $\exists$ over the disjunction and in each disjunct attach it to the cluster that its variable belongs to. This will go through because $\exists t\left[P(t)\land Q\right]\equiv\exists tP(t)\land Q$ when $Q$ does not depend on $t$. After renaming the expressions with  $\exists t$ attached into new predicates we will have a matrix in PDF and one less quantifier in the prefix. 

If the next quantifier is $\forall$ we will transform the matrix into PCF as in \eqref{distrib}, and push $\forall$ in the same way. If $\forall$ was the innermost one we will start from CNF instead. We repeat the process with each quantifier in the prefix. By induction on their number, at the end of it we will get an equivalent PDF over the initial partition.}. 

This opens a path to semantic description of disconnected relations, analogous to the one for degenerate relations in Section \ref{meas}. They are conglomerates of finitely many degenerate relations with the same relata independent of each other in each. The possibilities represented by the disjuncts need not be mutually exclusive, although one can always modify the PDF to make it so, at the expense of exploding the number of possibilities.

From PDFs it is easy to see that {\it on infinite domains} teridentity is a connected relation, and hence irreducible. Indeed, there are infinitely many domain elements, but only finitely many disjuncts in \eqref{oranddec}. By the pigeonhole principle, at least one disjunct will have to be true on infinitely many triples. But if $R$ is teridentity then restricting it to the subset of the domain with elements from those triples will give us teridentity on that subset. And it will be represented by $P_i(x,y)\land Q_i(z)$ for some $i$ there, i.e. by a Cartesian product. That is impossible, as we already know. The same goes for any other partition of places, and the same argument applies to all pluridentities. Thus, on infinite domains all pluridentities are connected. 

This completes our quest for an invariant form of PRT, at least on infinite domains. Every higher polyad reduces to triads by bonding alone, and no connected triad reduces to monads and dyads by any logical operations once they are explicated. 

Our excitement may be cooled off a notch once we realize that {\it any finite relation}, in particular, any relation on a finite domain, is disconnected, including teridentity and higher pluridentities. This is because we can dissociate finite relations into individual triples, define different $P_i$ and $Q_i$ for each, and then take the disjunction over them all. Specifically, using the delta-function monads $\delta_a(x):=I_2(a,x)$ we can decompose any finite relation as
\beq\label{deltadec}
R(x_1,\dots,x_n)\equiv\bigvee_{\vec{a}\,\in R}
\delta_{\vec{a}^{1}}(x_1)\land\dots\land\delta_{\vec{a}^{n}}(x_n),
\eeq
where $\vec{a}^{\,i}$ is the $i$-th entry of the tuple $\vec{a}$. 

One may notice a resemblance between the delta-function decomposition \eqref{deltadec} and the hypostatic abstraction \eqref{hypred}, which replaces the index $\vec{a}$ with a new variable $t$ that runs over tuples and the disjunction with $\exists\,t$. The difference is that we do not need to add new elements to the domain, or code them into the domain elements when defining $R^{(i)}(t,x)$. The latter can only be done for small relations anyway. In contrast, there is no such restriction on the number of disjuncts in \eqref{deltadec}, it only has to be finite. This means that delta-function decompositions remove the restriction to small relations in the reducibility clause of PRT.
\begin{figure}[!ht]
\begin{centering}
a)\ \ \ \includegraphics[width=0.25\textwidth]{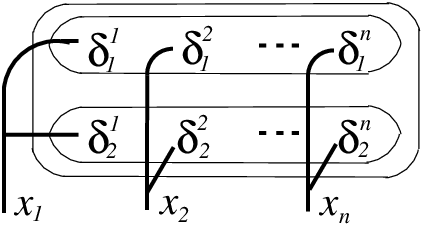} \hspace{0.2\textwidth} b)\ \ \ \includegraphics[width=0.15\textwidth]{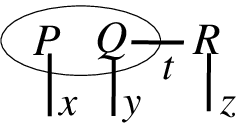}  \par
\end{centering}
\caption{\label{Delta} a) EG of the delta-function decomposition of an $n$-ad with two tuples. All triple points must be outside of cuts, we abbreviated $\delta_{\vec{a}^{\,i}_{\!j}}$ as $\delta^{i}_{\!j}$; b) EG template for a reducible triad that is neither additively nor multiplicatively degenerate.}
\end{figure}
\begin{quote} {\bf General reducibility (invariant form):} Any relation of adicity three or more is reducible to triads by explicated logical operations with $\exists,\land,\lor$.
\end{quote}
The cost is a massive increase in the number of required triads, see Figure \ref{Delta}\,a). PDFs generally, when explicated, require a multitude of them to identify variables in different disjuncts. 
As a result, delta-function decompositions cannot reduce teridentity even on finite domains, nor can they represent $n$-identities with fewer than $n-2$ triads. Some disconnected triads may still be irreducible. As a matter of fact, \cite{CorPos06,CorPos11} proved that teridentities are  generally irreducible even on finite domains, but the proof is too technical to discuss here.

\section{Which triads are genuine?}\label{gen}

What is striking about the general irreducibility clause of PRT is not just the technical complexity of its proof, but also that it is more than elaboration of the core argument that bringing two together takes a third. In non-positive logic it {\it does not} take a third to bring two (or even three) together, at least, not in the sense of using triads. Cuts can emulate them to a degree. Cartesian sums, that looked like genuine triads in the positive fragment, can be shattered into monads and held together by a cut, Figure \ref{EGcut}\,b). The simple valency reasons that prompted Peirce to call the irreducibility clause a ``truism" (\citetalias{CP}\,5.469, c.\,1906) do not suffice to explain it in general. Something else is at work in addition to them. 

Removing the restriction to positive operations in Definition \ref{ter+}, we obtain the definition of (general) 
{\it ternarity}, $\ter$, of a relation as the minimal number of triads in its explicated representation. Since the general PRT is stronger it may seem that ternarity supersedes everything positive ternarity had to say, but this is not the case; $\ter$ and $\ter^+$ are different relation invariants, and tell us different things about what goes on inside relations. We saw that in positive EG triads mediated information sharing between different relata, but in EG with cuts their role is opaquer. In addition to ``welding" different relata, triads are also used to ``weld" inputs of the {\it same} relatum under different possibilities covered by the relation. This is especially conspicuous in delta-function decompositions, see Figure \ref{Delta} a), although they are almost never reductions of  minimal ternarity. 

Minimal positive reductions are given by hypostatic abstraction for small relations , and for connected relations they are even generally minimal (by the same topological argument as before), but we do not even know if finite $n$-identities can be represented with fewer than $n-2$ triads when cuts are allowed. We do not yet have a clear picture of minimal general reductions to confidently say what they tell us about the role of triads, so we will only venture a preliminary conjecture. Some of the seemingly mediating triads in positive reductions convert into single relatum triads in general reductions, and those can be eliminated. When this happens $\ter$ is strictly smaller than $\ter^+$. This does not happen in connected relations, all of their triads are pure mediators, and it happens most conspicuously in monadic Cartesian sums. In a way, triads are traded for cuts. 

Peirce sometimes used ``irreducible" and ``non-degenerate" interchangeably when describing ``genuine" triads. Burch set up a system where these two notions come out as equivalent, but, as we saw, it comes at the price of severely weakening PRT. The two must come apart, and it is clear that irreducibility takes precedence. Thus, triads are {\it genuine} whenever they are irreducible. But which triads are irreducible, what {\it makes} a triad genuine? We are lacking a semantic description here, something discernible from relation's extension rather than from algebra\pagenote{One may object that extensional semantics is at odds with Peirce's conceptual view of relations. However, even intensional analysis, such as Burch's, needs it as a stepping stone. Algebraic properties are even more conceptually opaque.}.

Non-degeneracy and connectedness are semantic approximations of genuineness from above and below, respectively. Degenerate triads are all reducible even positively, and hence non-genuine, but disconnected triads may be genuine, like finite teridentities. \cite{CorPos06,CorPos11} discovered a finer property of non-genuine disconnected triads: of two monads with the same free variable in different disjuncts of their PDF one is always a subset of the other. This allowed them to prove genuineness of finite teridentities, but it is far from an  exhaustive description. It would be of interest to know how other semantic properties of relations, such as existence of a relatum that uniquely determines the rest (a {\it key}), interact with genuineness. In teridentities every relatum is a key.

We also do not know if infinite genuine triads can be disconnected, but we do know that non-genuine triads can fake genuineness pretty well. For example, there are many non-genuine triads that are neither additively nor multiplicatively degenerate. Their reductions, while disconnected, tangle bonds and cuts in a way that prevents pulling apart components with different relata.  The simplest candidate is shown on Figure \ref{Delta}\,b), and with `generic' enough $P$, $Q$, $R$ one can make it work\pagenote{In FOL the formula is $\exists t\left[\neg\left(P(x)\land Q(y,t)\right)\land R(t,z)\right]$. Taking $P(x):=(x\geq0)$, $Q:=\neg I_2$ and $R(t,z):=(t^2+z^2\leq1)$ we can express the relation more transparently (but implicitly using teridentity to identify $z$-s):
$$(x\geq0)\land(z^2\leq1)\lor(y^2+z^2\leq1).$$
The reader is invited to check that neither it nor its negation are Cartesian products.}.

Large (more tuples than domain elements) finite relations present a special challenge in Peircean analysis because the usual hypostatic abstraction does not work for them. Delta-function decompositions establish their reducibility when $\lor$ can be used in addition to $\exists,\land$, but the question remains whether they are {\it positively} reducible. \cite{Herz} gave an example of a non-degenerate tetrad on a $3$-element domain that is not a bond of two triads, as any small tetrad is by hypostatic abstraction. However, this does not rule out more complex reductions, and, indeed, one can show that Herzberger's tetrad is a bond of four triads. While this is more of an issue for higher polyads, because triads can always be `represented' by their single predicate letter, even for large triads, such as $\neg I_3$, we do not know if they are positively reducible to teridentities. 

And it is natural to extend the notion of genuineness to higher polyads. Of course, we cannot ask that they be irreducible, but we can ask that their reductions involve no less than $n-2$ triads, the maximum for small $n$-ads. 
It is open even whether finite $n$-identities for $n>3$ are genuine. Their ternarity is at least $1$, because $I_3$ can be bonded from them and it is known to be genuine, and it is at most $n-2$, because they are small. But is it exactly $n-2$ or can they be represented with fewer triads?

If a loose end is disconnected from the rest in a positive EG we can consolidate the entire component that contributes to that end into a new monad, but cuts partition the sheet of assertion in a way of their own and may block such consolidation. As a result, cuts considerably complicate determination of ternarity, and this, in turn, complicates philosophical analysis of Thirdness in non-positive logic.

\section{Conclusions}

With the general invariant form of PRT we have the negative answer to the title question -- PRT is not gerrymandered, it is a deep fact about the fine structure of relations that can be neither manufactured nor undone by manipulating algebraic conventions. 

The opposite impression was sustained by technical assumptions in the early proofs that were later removed. Such assumptions are only natural in initial proofs of a subtle conjecture, which PRT turned out to be. Unfortunately, some of them, such as the valency rule or the ban on negating juxtapositions, were taken for essential features, resulting in weaker and philosophically vulnerable versions of PRT that only reinforced the gerrymandering objection. In particular, we argued that  privileging some relational operations and/or graphical calculus, although traceable to some of Peirce's texts, does not fit into his broader epistemological framework. Our approach was instead to explicate relational constructions rather than to restrict them.

Our treatment clarified the role of EG and PAL as privileged {\it expository} devices for employment of triads in relational constructions rather than as ``fundamentally superior" to FOL, as Skidmore ascribed to Peirce, or as two among many ``constructive resources", with PRT becoming a syntactic artifact of favored formalisms. On the other hand, their expository privilege is easy enough to defend, just as the like privilege of eigenbases in analysis of linear operators. We submit that where Peirce did argue superiority of EG (e.g. in \citetalias{CP}\,4.561) it is this kind of superiority that he had in mind. 

In the traditional treatments of PRT, the focus was, essentially, on positive reductions, negation was taken along for the ride only to the extent that it did not get in the way. The removal of technical restrictions exposed new conceptual phenomena when non-positive operations (disjunction and negation) are used, that Peirce, and even Herzberger and Burch, did not get a chance to consider. Neither the valency rule nor hypostatic abstraction suffice for general reductions, and it is not so easy to trace the proofs of general PRT to the original intuition that bringing together two requires a third. If we are right that the full strength of general PRT is called for in Peirce's development of the categories then the philosophical import of general reductions needs to be investigated thoroughly.

We only took a first bite at it by introducing ternarity, a relation invariant inherent in PRT, and formulating new questions about genuineness of triads and polyads in the light of it. The chasm between reducibility properties of infinite and finite relations also warrants further investigation. Although all pluridentities are likely genuine, the question is still open on finite domains! It is also open whether some large finite relations are positively reducible to triads at all. 
Upon the explication, the semiotic role of triads as mediators of information sharing between different relata in positive reductions emerged very clearly. Understanding their role in general reductions, especially in connection with the dependencies studied in the database theory, is another promising track. Furthermore, as far as we know, all discussions of PRT so far have been restricted to first order operations. Higher order operations (such as transitive closures) are also actively studied, and naturally come up as database queries \citep{Buss01}. In what form and to what extent does PRT extend to them? The intriguing perspectives on semiotic analysis of relations and their Peircean classification that the invariant PRT points to present a rich field for further study.

{\footnotesize
\printnotes
}


{\small
\bibliographystyle{plainnat}
\bibliography{PeirceRed}
}

\end{document}